\documentclass{article}%
\usepackage{amsmath}
\usepackage{amsfonts}
\usepackage{amssymb}
\usepackage{graphicx}%
\providecommand{\U}[1]{\protect\rule{.1in}{.1in}}
\begin{document}

\title{Statistical Reasoning: Choosing and Checking the Ingredients, Inferences Based
on a Measure of Statistical Evidence with Some Applications}
\author{Luai Al-Labadi\thanks{Department of Mathematics, University of Sharjah},
Zeynep Baskurt\thanks{Genetics and Genome Biology, Hospital for Sick Children}
and Michael Evans\thanks{Department of Statistical Sciences, University of
Toronto}}
\date{}
\maketitle

\noindent\textit{Abstract}: The features of a logically sound approach to a
theory of statistical reasoning are discussed. A particular approach that
satisfies these criteria is reviewed. This is seen to involve selection of a
model, model checking, elicitation of a prior, checking the prior for bias,
checking for prior-data conflict and estimation and hypothesis assessment
inferences based on a measure of evidence. A long-standing anomalous example
is resolved by this approach to inference and an application is made to a
practical problem of considerable importance which, among other novel aspects
of the analysis, involves the development of a relevant elicitation
algorithm.\medskip

\bigskip

\section{Introduction}

It is relevant to ask what characteristics should be required of a theory of
statistical reasoning. The phrase \textit{statistical reasoning} is used here,
as opposed to statistical inference, because there is a logical separation
between how the ingredients to a statistical problem are chosen and checked
for their validity, and the inference step which involves the application of
the rules of a theory of inference to the ingredients. So it is argued in
Section 2 that there are two aspects to a theory of statistical reasoning: (i)
specifying methodology for choosing and checking the ingredients to a
statistical analysis beyond the data and (ii) specifying a theory of inference
using these ingredients that is based on a measure of statistical evidence.
These components correspond to the premises and the argument in logical reasoning.

In Section 3 a specific theory of statistical reasoning, as described in Evans
(2015), that satisfies the criteria developed in Section 2 is reviewed. It is
shown that an application of the theory of relative belief inference resolves
difficulties in a problem that has led to anomalous results for other
theories. The overall aim of the paper is to argue in favor of this approach
to statistical reasoning based on its logical coherence and its utility in applications.

In Section 4 the theory of statistical reasoning is applied to an important
practical problem where some inferential difficulties have arisen, namely, the
problem of determining whether or not there is a relationship between a
binary-valued response variable and predictors. For this, response
$y\in\{0,1\}$ is related to $k$ quantitative predictors $\mathbf{x}%
=(x_{1},\ldots,x_{k})$ via $y\sim$ Bernoulli$(p(\mathbf{x}))$ with
\begin{equation}
p(\mathbf{x})=G(\beta_{1}x_{1}+\cdots+\beta_{k}x_{k}) \label{eq1}%
\end{equation}
where $G$ is a known cdf and $\mathbf{\beta}=(\beta_{1},\ldots,$ $\beta
_{k})\in R^{k}$ is unknown. This can be regarded as a case-study for the
overall approach and a number of novel results are obtained. Perhaps the
biggest challenge with this model is determining a suitable prior and in
Section 4.1 an elicitation algorithm is developed that improves on earlier
efforts. The bias in the prior is measured in Section 4.2. Model checking is
essential, namely, determining if (\ref{eq1}) holds at least approximately.
Since this is dealt with in Al-Labadi, Baskurt and Evans (2017), that approach
is used here without much discussion. The check for prior-data conflict is
developed in Section 4.3, together with an approach to modifying the prior
when necessary, and relative belief inferences are applied in Section 4.4.

\section{The Foundations of Statistical Reasoning}

When concerned with the foundations of statistics it is reasonable to ask:
what is the purpose of statistics as a subject or what is its role in science?
To answer this, consider a context where an investigator has interest in some
quantity and either wants to know (E) the value of this quantity or has a
theory that leads to a specific value for the quantity and so wants to know
(H)\ if this value is indeed correct and so test the theory. In addition, the
investigator has available data $d,$ produced in some fashion, which it is
believed contains \textit{evidence} concerning answers to (E) and (H). The
purpose of statistical theory is to provide a reasoning process that can be
applied to $d$ to determine what the evidence has to say about (E) or (H),
namely, produce an estimate of the quantity based on the evidence or assess
whether there is evidence either in favor of or against the hypothesized
value. Also, as is widely recognized, estimation and hypothesis assessment
should also produce a measure of the accuracy of the estimate and a measure of
the strength of the evidence for or against the hypothesis. Answering (E)
and/or (H) is called statistical inference and a sound, logical theory of
statistical inference, that contained the minimal ingredients possible, can be
viewed as a major goal of the subject.

Any theory that does not lead to specific answers to (E) and (H) or is
dependent on ingredients or rules of reasoning that are not well-justified, is
unsatisfactory. In the end, the believability of the inferences drawn is
entirely dependent on the soundness of the theory which produced them. So
statistics is not an empirical subject like physics, where conclusions can
also be assessed against the empirical world, but is more like an extension of
purely logical reasoning to contexts where the data does not lead to
categorical answers to (E) and (H) and so produces uncertainty. The view is
taken here that we want to maintain a close relationship between a theory of
statistical reasoning and the theory of logical reasoning. This has a number
of consequences with perhaps the most important being that it implies a
separation of the assessment of the appropriateness of the ingredients
specified to a statistical analysis beyond the data, and the theory producing
the inferences. The ingredients play the role of the premises and the theory
of statistical inference takes the role of the rules of inference used in a
logical argument. The separation of these aspects of a logical argument has
been understood since Aristotle, see Kneale and Kneale (1962).

There are two main theses of the argument developed here (i) all ingredients
to a statistical analysis must be checkable against the data and (ii) the
theory of inference must be based on a measure of statistical evidence. The
rationale for (i) and (ii) are now considered with (ii) discussed first, as it
plays a key role.

The concept of the evidence in the data is clearly of utmost importance to
statistical reasoning. There is no need, however, to provide a measure of the
total evidence contained in the data. For the measure of evidence only has to
deal with (E) and (H) for the quantity of interest. The measure of evidence
must clearly indicate whether there is evidence for or against any specific
value of the quantity of interest being the true value. This follows from\ the
desired relationship with logic, as the rules of logical inference assume the
truth of the premises, so the theory of statistical inference has to be based
on the truth of the ingredients and this implies that one of the possible
values for the quantity of interest is true. A theory of logical reasoning
that could only determine falsity and never truth is not useful and similarly
any valid measure of evidence must be able to indicate evidence in favor as
well evidence against.

Once a measure of statistical evidence is determined, an estimate of the
quantity of interest is necessarily the value that maximizes the measure of
evidence and the accuracy of the estimate can be assessed by looking at the
size of the set of values that have evidence in their favor. The measure of
evidence similarly necessarily determines whether there is evidence for or
against a hypothesized value and the strength of this evidence can be assessed
by comparison with the evidence associated with the other possible values for
the quantity of interest.

Consider now requirement (i). If a satisfactory measure of evidence could be
determined from the data alone, then this would be ideal, but currently this
is not available and it is questionable whether it is even possible. It is
\textit{assumed} hereafter that the data $x\in\mathcal{X}$ can be regarded as
having been produced by some probability distribution on the set $\mathcal{X}$
with unknown density $f.$ If the data was collected via random sampling, then
this assumption seems justified, but it is always an assumption. The density
$f$ is unknown and it assumed that once $f$ is known, then this determines
answers to (E) and ((H). The ingredients are then as follows: it is
\textit{assumed} that $f\in\{f_{\theta}:\theta\in\Theta\},$ a collection of
densities on $\mathcal{X}$ indexed by the parameter $\theta\in\Theta$ called
the statistical model, and it is \textit{assumed} that there is a probability
measure $\Pi$ with density $\pi$ on $\Theta$ that represents beliefs of the
investigator about the true value of $\theta\in\Theta$ and called the prior.
The ingredients correspond to the premises of a logical argument and these may
be true or false.

It can be questioned as to whether both the model and prior are necessary for
the development of a satisfactory theory and certainly minimizing the
ingredients is desirable. But as discussed in Section 3 it seems that a valid
definition of a measure of evidence requires both and again the challenge is
open to develop a satisfactory measure of evidence that uses fewer
ingredients. In particular, the use of a prior is often claimed to be
inappropriate as it is subjective in nature and, as the goal of a scientific
investigation is to be as objective as possible, the prior seems contrary to
that. It needs to be recognized, however, that all ingredients to a
statistical analysis beyond the data are subjective in nature as they are
chosen by the statistician. As discussed in Section 3.2, it is possible to
check both the model and the prior against the (objective) data to determine
whether or not reasonable choices have been made. Also, it is possible to
check whether or not the chosen ingredients have biased the results so that
the inferences obtained are in fact foregone conclusions, namely, could have
been made without even looking at the data. It is our view that checking for
bias and checking for conflict with the data go a long way towards answering
criticisms concerning the subjectivity inherent in a statistical analysis.
Another implication from (i) is that no ingredient can be added to a
statistical analysis unless it can be checked against the data which rules out
the use of loss functions.

It is not clear how the ingredients are to be chosen and guidance needs to be
provided for this. Not much has been written about how the model is to be
chosen but certainly something needs to be said to justify a specific choice
as part of the statistical reasoning argument. Much more has been written
about the selection of the prior and the position is adopted here that it is
necessary to base this on a clearly stated elicitation algorithm, namely, a
prescription for how an expert can translate knowledge into beliefs as
expressed via the prior.

In summary, the desiderata for a theory of statistical reasoning include the
following: a methodology for choosing a model, an elicitation algorithm for
selecting a prior, methodology for assessing the bias in the ingredients
chosen, model checking and checking for prior-data conflict procedures and a
theory of inference based upon a measure of statistical evidence.

\section{A Theory of Statistical Reasoning}

Choosing and checking the ingredients logically comes before inference but it
is convenient to discuss these in reverse order.

\subsection{Relative Belief Inferences}

Consider now a statistical problem with ingredients the data $d,$ a model
$\{f_{\theta}:\theta\in\Theta\},\ $a prior $\pi$ and interest is in making
inference about $\psi=\Psi(\theta)$ for $\Psi:\Theta\rightarrow\Psi$ where no
distinction is made between the function and its range to save notation.
Initially, suppose that all the probability distributions are discrete. This
isn't really a restriction in the discussion, however, as if something works
for inference in the discrete case but does not work in the continuous case,
then it is our view that the concept is not being applied correctly or the
mathematical context is just too general. For us the continuous case is always
to be thought of as an approximation to a fundamentally discrete context, as
measurements are always made to finite accuracy, and the approximation arises
via taking limits. Some additional comments on the continuous case are made subsequently.

As discussed in Section 2, the basic object of inference is the measure of
evidence and what is wanted is a measure of the evidence that any particular
value $\psi\in\Psi$ is true. Based on the ingredients specified, there are two
probabilities associated with this value, namely, the prior probability
$\pi_{\Psi}(\psi)$, as given here by the marginal prior density evaluated at
$\psi,$ and the posterior probability $\pi_{\Psi}(\psi\,|\,d)$, as given here
by the marginal posterior density evaluated at $\psi.$ In certain treatments
of inference $\pi_{\Psi}(\psi\,|\,d)$ is taken as a measure of the evidence
that $\psi$ is the true value but, for a wide variety of reasons, this is not
felt to be correct and Example 1 provides a specific case where this fails.
Also, this measure suffers from the same basic problem of $p$-values, namely,
there is no obvious dividing line between evidence for and evidence against.
Moreover, it is to be noted that probabilities measure belief and not
evidence. If we start with a large prior belief in $\psi$, then unless there
is a large amount of data, there will still be a large posterior belief even
if it is false and similarly, if we started with a small amount of belief.
There is agreement, however, to use the \textit{principle of conditional
probability} to update beliefs after receiving information or data and this is
to be regarded as the first principle of the theory of relative belief.

So how is the evidence that $\psi$ is the true value to be measured? Basic to
this is the \textit{principle of evidence}: if $\pi_{\Psi}(\psi\,|\,d)>\pi
_{\Psi}(\psi)$ there is evidence in favor as belief has increased, if
$\pi_{\Psi}(\psi\,|\,d)<\pi_{\Psi}(\psi)$ there is evidence against as belief
has decreased and if $\pi_{\Psi}(\psi\,|\,d)=\pi_{\Psi}(\psi)$ there is no
evidence either way. This principle has a long history in the philosophical
literature concerning evidence and a nice discussion can be found in Salmon
(1973). This principle doesn't provide a specific measure of evidence but at
least it indicates clearly when there is evidence for or against, independent
of the size of initial beliefs, and it does suggest that any reasonable
measure of the evidence depends on the difference, in some sense, between
$\pi_{\Psi}(\psi)$ and $\pi_{\Psi}(\psi\,|\,d),$ namely, evidence is measured
by \textit{change in belief} rather than belief. A number of measures of this
change have been proposed, see Evans (2015) for a discussion, but by far the
simplest and the one that has the nicest properties is the relative belief
ratio%
\begin{equation}
RB_{\Psi}(\psi\,|\,d)=\pi_{\Psi}(\psi\,|\,d)/\pi_{\Psi}(\psi). \label{eqRB}%
\end{equation}
So if $RB_{\Psi}(\psi\,|\,d)>(<,=)1$ there is evidence for (against, neither)
for $\psi$ being the true value.\ The use of the relative belief ratio to
measure the evidence is the third and final principle of the theory, what we
call the \textit{principle of relative belief}. The relative belief ratio can
also be written as $RB_{\Psi}(\psi\,|\,d)=m(d\,|\,\psi)/m(d)$ where $m$ is the
prior predictive density of the data and $m(\cdot\,|\,\psi)$ is the
conditional prior predictive density of the data given $\Psi(\theta)=\psi.$

Another natural candidate for a measure of evidence is the Bayes factor
$BF_{\Psi}(\psi\,|\,d)$ as this satisfies the principle of evidence, namely,
$BF(\psi\,|\,d)>(<,=)1$ when there is evidence for (against, neither) $\psi$
being the true value. The Bayes factor can be defined in terms of the relative
belief ratio as $BF_{\Psi}(\psi\,|\,d)=RB_{\Psi}(\psi\,|\,d)/RB_{\Psi}%
(\{\psi\}^{c}\,|\,d)$ but not conversely. Note that the relative belief ratio
of a set $A\subset\Psi$ is $RB_{\Psi}(A\,|\,d)=\Pi_{\Psi}(A\,|\,d)/\Pi_{\Psi
}(A)$ where $\Pi_{\Psi},\Pi_{\Psi}(\cdot\,|\,d)$ are the prior and posterior
probability measures of $\Psi,$ respectively. It might appear that $BF_{\Psi
}(\psi\,|\,d)$ is a comparison between the evidence for $\psi$ being true with
the evidence for $\psi$ being false but it is provable that $RB_{\Psi
}(A\,|\,d)>1$ implies $RB_{\Psi}(A^{c}\,|\,d)<1$ and conversely, so this is
not the case. Also, as subsequently discussed, in the continuous case it is
natural to take $BF_{\Psi}(\psi\,|\,d)=RB_{\Psi}(\psi\,|\,d).$

The principle of relative belief leads to an ordering of the possible values
for $\psi$ as $\psi_{1}$ is preferred to $\psi_{2}$ whenever $RB_{\Psi}%
(\psi_{1}\,|\,d)>RB_{\Psi}(\psi_{2}\,|\,d)$ since there is more evidence for
$\psi_{1}$ than $\psi_{2}.$ When $\Psi(\theta)=\theta$ this agrees with the
likelihood ordering but likelihood fails to provide such an ordering for
general $\psi.$ It is common to use the profile likelihood ordering even
though this is not a likelihood ordering and this doesn't agree with the
relative belief ordering. It is noteworthy that the relative belief idea is
consistent in the sense that inferences for all $\psi=\Psi(\theta)$ are based
on a measure of the change in prior to posterior probabilities.

The relative belief ordering leads immediately to a theory of estimation. For
basing inferences on the evidence requires that the relative belief estimate
be a value $\psi(d)$ maximizes $RB_{\Psi}(\psi\,|\,d)$ and typically this
value is unique so $\psi(d)=\arg\sup_{\psi\in\Psi}RB_{\Psi}(\psi\,|\,d)$. It
is also necessary to say something about the accuracy of this estimate in an
application. For this a set of values containing $\psi(d)$ is quoted and the
\textquotedblleft size\textquotedblright\ of the set is taken as the measure
of accuracy. Again following the ordering based on the evidence, it is
necessary that the set take the form $\{\psi:RB_{\Psi}(\psi\,|\,d)>c\}$ for
some constant $c\leq\sup_{\psi\in\Psi}RB_{\Psi}(\psi\,|\,d)$ since, if
$RB_{\Psi}(\psi_{1}\,|\,d)\leq RB_{\Psi}(\psi_{2}\,|\,d),$ then $\psi_{2}$
must be included whenever $\psi_{1}$ is. But what $c$ should be used? It is
perhaps natural to chose $c$ so that $\{\psi:RB_{\Psi}(\psi\,|\,d)>c\}$
contains some prescribed amount of posterior probability, so the set is a
$\gamma$-credible region. But there are several problems with this
approach.\ For what $\gamma$ should be chosen? Even if one is content with
some particular $\gamma,$ say $\gamma=0.95,$ there is the problem that the set
may contain values $\psi$ with $RB_{\Psi}(\psi\,|\,d)<1$ and such a value has
been ruled out since there is evidence against such a $\psi$ being true. It is
argued in Evans (2015) that the \textit{plausibility set} $Pl_{\Psi}%
(d)=\{\psi:RB_{\Psi}(\psi\,|\,d)>1\}$ be used as $Pl_{\Psi}(d)$ contains all
the values for which there is evidence in favor of it being the true value. In
general circumstances, it is provable that $RB_{\Psi}(\psi(d)\,|\,d)>1$ so
$Pl_{\Psi}(x)\neq\phi.$ There are several possible measures of size but
certainly the posterior content $\Pi_{\Psi}(Pl_{\Psi}(d)\,|\,d)$ is one as
this measures the belief that the true value is in $Pl_{\Psi}(d),$ but also
some measure such as length or cardinality is relevant. If $Pl_{\Psi}(d)$ is
small and $\Pi_{\Psi}(Pl_{\Psi}(d)\,|\,d)$ large, then $\psi(d)$ can be judged
to be an accurate estimate of $\psi.$

It is immediate that $RB_{\Psi}(\psi_{0}\,|\,d)$ is the evidence concerning
$H_{0}:\Psi(\theta)=\psi_{0}.$ The evidential ordering implies that the
smaller $RB_{\Psi}(\psi_{0}\,|\,d)$ is than 1, the stronger is the evidence
against $H_{0}$ and the bigger it is than 1, the stronger is the evidence in
favor $H_{0}.$ But how is one to measure this strength? In Baskurt and Evans
(2013) it is proposed to measure the \textit{strength of the evidence} via%
\begin{equation}
\Pi_{\Psi}\left(  \left.  RB_{\Psi}(\psi\,|\,d)\leq RB_{\Psi}(\psi
_{0}\,|\,d)\,\right\vert \,d\right)  \label{eq2}%
\end{equation}
which is the posterior probability that the true value of $\psi$ has evidence
no greater than that obtained for the hypothesized value $\psi_{0}.$ When
$RB_{\Psi}(\psi_{0}\,|\,d)<1$ and (\ref{eq2}) is small, then there is strong
evidence against $H_{0}$ since there is a large posterior probability that the
true value of $\psi$ has a larger relative belief ratio. Similarly, if
$RB_{\Psi}(\psi_{0}\,|\,d)>1$ and (\ref{eq2}) is large, then there is strong
evidence that the true value of $\psi$ is given by $\psi_{0}$ since there is a
large posterior probability that the true value is in $\{\psi:RB_{\Psi}%
(\psi\,|\,x)\leq RB_{\Psi}(\psi_{0}\,|\,d)\}$ and $\psi_{0}$ maximizes the
evidence in this set. Additional results concerning $RB_{\Psi}(\psi
_{0}\,|\,d)$ as a measure of evidence and (\ref{eq2}) can be found in Baskurt
and Evans (2013) and Evans (2015).

For continuous parameters it is natural to define $RB_{\Psi}(\psi
\,|\,d)=\lim_{\epsilon\rightarrow0}RB_{\Psi}($\newline$N_{\epsilon}%
(\psi)\,|\,d)$where $N_{\epsilon}(\psi)$ is a sequence of sets converging
nicely to $\{\psi\}$ as $\epsilon\rightarrow0.$ When the densities are
continuous at $\psi,$ then this limit equals (\ref{eqRB}) so this is a measure
of evidence in general circumstances. Also, it is natural to define the Bayes
factor by $BF_{\Psi}(\psi\,|\,d)=\lim_{\epsilon\rightarrow0}BF_{\Psi
}(N_{\epsilon}(\psi)\,|\,d)$ and, when the densities are continuous at $\psi,$
then $BF_{\Psi}(\psi\,|\,d)=RB_{\Psi}(\psi\,|\,d).$

A variety of consistency results, as the amount of data increases, are proved
in Evans (2015) concerning the estimation and hypothesis assessment
procedures. In particular, when $H_{0}$ is false, then (\ref{eqRB}) converges
to 0 as does (\ref{eq2}). When $H_{0}$ is true, then (\ref{eqRB}) converges to
its largest possible value (greater than 1 and often $\infty$) and, in the
discrete case (\ref{eq2}) converges to 1. In the continuous case, however,
when $H_{0}$ is true, then (\ref{eq2}) typically converges to a $U(0,1)$
random variable. This simply reflects the approximate nature of the inferences
and is easily resolved by requiring that a deviation $\delta>0$ be specified
such that if dist$(\psi_{1},\psi_{2})<\delta,$ where dist is some measure of
distance determined by the application, then this difference is to be regarded
as immaterial. This leads to redefining $H_{0}$ as $H_{0}=\{\psi:$
dist$(\psi,\psi_{0})<\delta\}$ and typically a natural discretization of
$\Psi$ exists with\ $H_{0}$ as one of its elements. With this modification
(\ref{eq2}) converges to 1 as the amount of data increases when $H_{0}$ is
true. Given that data is always measured to finite accuracy, the value of a
typical continuous-valued parameter can only be known to a certain finite
accuracy no matter how much data is collected. So such a $\delta$ always
exists and it is part of an application to determine the relevant value, see
Example 7 here, Al-Labadi, Baskurt and Evans (2017) and Evans, Guttman and Li
(2017) for developments on determining $\delta$.

It is immediate that relative belief inferences have some excellent
properties. For example, any 1-1 increasing function of $RB_{\Psi}%
(\cdot\,|\,d),$ such as $\log RB_{\Psi}(\cdot\,|\,d),$ can be used to measure
evidence as the inferences are invariant to this choice. Also, $RB_{\Psi
}(\cdot\,|\,d)$ is invariant under smooth reparameterizations and so all
relative belief inferences possess this invariance property. For example,
MAP\ (maximum a posteriori)\ inferences are not invariant and this leads to
considerable doubt about their validity, see also Example 1. In Evans (2015)
results from a number of papers are summarized establishing optimality results
for relative belief inferences in the collection of all Bayesian inferences.
For example, Al-Labadi and Evans (2017) establish that relative belief
inferences for $\psi$ have optimal robustness to the prior $\pi_{\Psi}$
properties. Also, as discussed in Section 3.2, since the inferences are based
on a measure of evidence a key criticism of Bayesian methodology can be
addressed, namely, the extent to which the inferences are biased can be measured.

Relative belief prediction inferences for future data are naturally produced
by using the ratio of the posterior to prior predictive densities for the
quantity in question. The following example illustrates this and demonstrates
significant advantages for relative belief.\smallskip

\noindent\textbf{Example 1.}\textit{ Prediction for Bernoulli sampling.}

Consider an example discussed in\ Chapter 6 of Diaconis and Skyrms (2018). A
tack is flipped with $x=1$ indicating the tack finishes point up and $x=0$
otherwise, so $x\sim$ Bernoulli$(\theta).$ Suppose the prior is $\theta\sim
U(0,1)$ and the goal is to predict $f$ future observations $(y_{1}%
,\ldots,y_{f})$ having observed $n$ independent tosses $(x_{1},\ldots,x_{n})$.
The posterior of $\theta$ is beta$(n\bar{x}+1,n(1-\bar{x})+1),$ the prior
predictive density of $(x_{1},\ldots,x_{n})$ is $m_{n}(x_{1},\ldots
,x_{n})=1/(n+1)\binom{n}{n\bar{x}}$ and the posterior predictive density for
$(y_{1},\ldots,y_{f})$ is%
\begin{equation}
m_{n,f}((y_{1},\ldots,y_{f})\,|\,(x_{1},\ldots,x_{n}))=\frac{(n+1)\binom
{n}{n\bar{x}}}{(n+f+1)\binom{n+f}{(n+f)\left[  \frac{n}{n+f}\bar{x}+\frac
{f}{n+f}\bar{y}\right]  }}, \label{postpred1}%
\end{equation}
which is constant for all $(y_{1},\ldots,y_{f})$ with the same value of
$\bar{y}.$ Maximizing (\ref{postpred1}) gives the MAP predictor of
$(y_{1},\ldots,y_{f}).$ If $n\bar{x}/(n+f)>1/2,$ then the maximum occurs at
$(y_{1},\ldots,y_{f})$ with $\bar{y}=1,$ namely, $(y_{1},\ldots,y_{f}%
)=(1,\ldots,1).$ If $n\bar{x}/(n+f)<1/2,$ then the maximum occurs at
$(y_{1},\ldots,y_{f})$ with $\bar{y}=0,$ namely, $(y_{1},\ldots,y_{f}%
)=(0,\ldots,0).$ If $n\bar{x}/(n+f)=1/2$, then a maximum occurs at both
$(y_{1},\ldots,y_{f})=(0,\ldots,0)$ and $(y_{1},\ldots,y_{f})=(1,\ldots,1)$.
So using MAP gives the absurd result that $(y_{1},\ldots,y_{f})$ is always
predicted to be either all 0's or all 1's. Clearly there is a problem here
with using MAP.

Now suppose $(x_{1},\ldots,x_{n})=(0,\ldots,0)$ so the prediction is all 0's
and
\[
m_{n,f}((y_{1},\ldots,y_{f})\,|\,(0,\ldots,0))=(n+1)/(n+f+1)\binom{n+f}%
{f\bar{y}}.
\]
For fixed $f,$ then $m_{n,f}((y_{1},\ldots,y_{f})\,|\,(0,\ldots,0))\rightarrow
0$ as $n\rightarrow\infty$ whenever $\bar{y}\neq0$ and converges to 1 when
$\bar{y}=0$. Diaconis and Skyrms (2018) note, however, that when $f=n$, then
$m_{n,n}((0,\ldots,0)\,|\,(0,\ldots,0))\rightarrow1/2$ as $n\rightarrow\infty$
and make the comment \textquotedblleft If this is an unwelcome surprise, then
perhaps the uniform prior is suspect.\textquotedblright\ They also refer to
some attempts to modify the prior to avoid this phenomenon which clearly
violates an essential component of the Bayesian approach. In our view there is
nothing wrong with the uniform prior, rather the problem lies with using
posterior probabilities implicitly as measures of evidence, both to determine
the predictor and to assess its reliability.

\begin{figure}
[ptb]
\begin{center}
\includegraphics[
height=2.4163in,
width=2.4163in
]%
{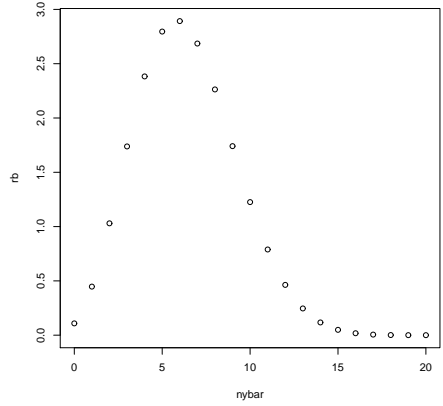}%
\caption{Plot of the relative belief ratio when $n=20,n\bar{x}=6$ in Example
1.}%
\label{diaconisfig}%
\end{center}
\end{figure}

The relative belief ratio for $(y_{1},\ldots,y_{f})$ is
\begin{equation}
RB((y_{1},\ldots,y_{f})\,|\,(x_{1},\ldots,x_{n}))=\frac{(n+1)\binom{n}%
{n\bar{x}}(f+1)\binom{f}{f\bar{y}}}{(n+f+1)\binom{n+f}{(n+f)\left[  \frac
{n}{n+f}\bar{x}+\frac{f}{n+f}\bar{y}\right]  }}. \label{relbel}%
\end{equation}
With $n=f=20$ and $n\bar{x}=6,$ Figure \ref{diaconisfig} gives the plot of
(\ref{relbel}) as a function of $n\bar{y}.$ The best relative belief predictor
of $(y_{1},\ldots,y_{f})$ is any sample with $f\bar{y}=6$ and $Pl_{n}%
(x_{1},\ldots,x_{n})=\{(y_{1},\ldots,y_{f}):f\bar{y}=2,3,\ldots,10\}$ has
posterior content $0.893.$ So there is reasonable belief that the plausibility
set contains the \textquotedblleft true\textquotedblright\ future sample but
certainly there are many such samples. By contrast with MAP,\ a sensible
prediction is made using relative belief.%
%TCIMACRO{\FRAME{ftbpFU}{2.4163in}{2.4163in}{0pt}{\Qcb{Plot of the relative
%belief ratio when $n=20,n\bar{x}=6$ in Example 1.}}{\Qlb{diaconisfig}%
%}{plot1.eps}{\special{ language "Scientific Word";  type "GRAPHIC";
%width 2.4163in;  height 2.4163in;  depth 0pt;  original-width 6.9998in;
%original-height 6.9998in;  cropleft "0";  croptop "1";  cropright "1";
%cropbottom "0";  filename 'plot1.eps';file-properties "XNPEU";}} }%
%BeginExpansion
%EndExpansion

For the case when $f=n$ and $(x_{1},\ldots,x_{n})=(0,\ldots,0),$ then%
\[
RB((y_{1},\ldots,y_{n})\,|\,(0,\ldots,0))=\frac{(n+1)^{2}}{(2n+1)}\frac
{\binom{n}{n\bar{y}}}{\binom{2n}{n\bar{y}}}=\frac{(n+1)^{2}}{(2n+1)}%
%TCIMACRO{\tprod _{i=0}^{n-1}}%
%BeginExpansion
{\textstyle\prod_{i=0}^{n-1}}
%EndExpansion
\left\{  \frac{2-\bar{y}-i/n}{(2-i/n)}\right\}
\]
which is decreasing in $\bar{y}$ and so is maximized for the sample with
$\bar{y}=0.$ Similarly, when $\bar{x}=1$ the predictor is the sample with
$\bar{y}=1.$ So, at the extremes, the predictions based on MAP\ and relative
belief are the same but otherwise there is a sharp disagreement. Also,
\[
Pl_{n}(0,\ldots,0)=\{(y_{1},\ldots,y_{n}):RB((y_{1},\ldots,y_{n}%
)\,|\,(0,\ldots,0))>1\}
\]
always contains $(y_{1},\ldots,y_{n})=(0,\ldots,0)$ and for any $c\in(0,1]$
such that $\bar{y}\geq c,$%
\begin{align*}
&  RB((y_{1},\ldots,y_{n})\,|\,(0,\ldots,0))=\{(n+1)^{2}/(2n+1)\}%
%TCIMACRO{\tprod _{j=0}^{n-1}}%
%BeginExpansion
{\textstyle\prod_{j=0}^{n-1}}
%EndExpansion
\left(  1-\bar{y}/(2-j/n)\right) \\
&  \leq\{(n+1)^{2}/(2n+1)\}\left(  1-c/2\right)  ^{n}\rightarrow0
\end{align*}
as $n\rightarrow\infty.$ Therefore, for any $c\in(0,1],$ there is an $N,$ such
that for all $n>N,$ then $Pl_{n}(0,\ldots,0)$ contains no $(y_{1},\ldots
,y_{n})$ having a proportion of 1's that is $c$ or greater. So $Pl_{n}%
(0,\ldots,0)$ is shrinking as $n$ increases in the sense that it contains only
samples with smaller and smaller proportion of 1's as $n$ increases.

The posterior content of the plausibility region equals%
\begin{equation}
\sum_{\{n\bar{y}:RB((y_{1},\ldots,y_{n})\,|\,(0,\ldots,0))>1\}}m_{n,f}%
((y_{1},\ldots,y_{f})\,|\,(0,\ldots,0))\binom{n}{n\bar{y}} \label{plausprob}%
\end{equation}
which equals the sum over all the summands that are greater than
$1/(n+1)\ $and%
\begin{align*}
&  m_{n,f}((y_{1},\ldots,y_{f})\,|\,(0,\ldots,0))\binom{n}{n\bar{y}}%
=\frac{(n+1)}{(2n+1)}%
%TCIMACRO{\tprod _{i=0}^{n-1}}%
%BeginExpansion
{\textstyle\prod_{i=0}^{n-1}}
%EndExpansion
\left\{  \frac{2-\bar{y}-i/n}{(2-i/n)}\right\} \\
&  =\frac{(n+1)}{(2n+1)}\left\{
\begin{array}
[c]{ccc}%
1 &  & \bar{y}=0\\
\frac{1}{2} &  & \bar{y}=\frac{1}{n}\\
\frac{1}{2}\frac{1-1/n}{2-1/n} &  & \bar{y}=\frac{2}{n}\\
\vdots &  & \vdots\\
\frac{1}{2}\frac{1-1/n}{2-1/n}\cdots\frac{1-(k-1)/n}{2-(k-1)/n} &  & \bar
{y}=\frac{k}{n}\\
\vdots &  & \vdots
\end{array}
\right.  .
\end{align*}
When $\bar{y}=k/n$ the corresponding term converges to $(1/2)^{k+1}.$ Thus for
all $n$ large enough, the sum (\ref{plausprob}) contains the terms for
$\bar{y}=0,1/n,\ldots,k/n.$ Therefore, for $\epsilon>0$ and all $n$ large
enough, (\ref{plausprob}) is greater than $(1/2)[1+1/2+\cdots+(1/2)^{k}%
]-\epsilon=1-(1/2)^{k+1}-\epsilon$ and the posterior content of $Pl_{n}%
(0,\ldots,0)$ converges to 1.

So relative belief also behaves appropriately when $f=n$ and $\bar{x}=0$ while
MAP does not. The failure of MAP might be attributed to the requirement that
the entire sample $(y_{1},\ldots,y_{n})$ be predicted. If instead it was
required only to predict the value $n\bar{y},$ then the prior predictive of
this quantity is uniform on $\{0,1,\ldots,f\},$ the posterior of $n\bar{y}$
equals $RB((y_{1},\ldots,y_{f})\,|\,(x_{1},\ldots,x_{n}))/(f+1)$ and the
relative belief ratio for $n\bar{y}$ equals $RB((y_{1},\ldots,y_{f}%
)\,|\,(x_{1},\ldots,x_{n})).$ So, as is often the case when the quantity in
question has a uniform prior, MAP and relative belief estimates are the same.
But even in this case there is no natural cut-off for MAP inferences to say
when there is evidence for or against a particular value. The fact that it is
necessary to modify the problem in this way to get a reasonable inference is,
in our view, a substantial failing of MAP. It seems reasonable to suggest that
when an inference approach is shown to perform poorly on such examples, that
it not be generally recommended. Additional examples of poor performance of
MAP are discussed in Evans (2015).

\subsection{Choosing and Checking the Ingredients}

The first choice that must be made is the model and there are a number of
standard models used in practice. There isn't a lot written about this step,
however, and yet it is perhaps the most important step in solving a
statistical problem. It is generally accepted that the correct way to choose a
prior is through elicitation. This means that a methodology is prescribed that
directs an expert in the application area on how to translate their knowledge
into a prior. There are various default priors in use that avoid this
elicitation step, but it is far better to recommend that sufficient time and
energy be allocated for the elicitation of a proper prior. Staying within the
context of probability suggests that a variety of paradoxes and illogicalities
are avoided.

Given the ingredients, the relative belief inferences may be applied correctly
but it is still reasonable to ask if these ingredients are appropriate for the
particular application. If not, then the inferences drawn cannot be considered
valid. There are at least two questions about the ingredients that need to be
answered, namely, is there bias inherent in the choice of ingredients and are
the ingredients contradicted by the data?

The concern for bias is best understood in terms of assessing the hypothesis
$H_{0}:\Psi(\theta)=\psi_{0}.$ Let $M(\cdot\,|\,\psi)$ denote the prior
predictive distribution of the data given that $\Psi(\theta)=\psi.$ Bias
against $H_{0}$ means that the ingredients are such that, with high
probability, evidence will not be obtained in favor of $H_{0}$ even when it is
true. Bias against is thus measured by
\begin{equation}
M(RB_{\Psi}(\psi_{0}\,|\,D)\leq1\,|\,\psi_{0}). \label{biasagainst}%
\end{equation}
If (\ref{biasagainst}) is large, then obtaining evidence against $H_{0}$ seems
like a foregone conclusion. For bias in favor of $H_{0}$ consider $M(RB_{\Psi
}(\psi_{0}\,|\,D)\geq1\,|\,\psi_{\ast})\ $where dist$(\psi_{\ast},\psi
_{0})=\delta,$ so $\psi_{\ast}$ is a value that just differs from the
hypothesized value by a meaningful amount. Bias in favor of $H_{0}$ is then
measured by
\begin{equation}
\sup_{\psi_{\ast}\in\{\psi:\text{ dist}(\psi,\psi_{0})=\delta\}}M(RB_{\Psi
}(\psi_{0}\,|\,D)\geq1\,|\,\psi_{\ast}). \label{biasfor}%
\end{equation}
If (\ref{biasfor}) is large, then obtaining evidence in favor of $H_{0}$ seems
like a foregone conclusion. Typically $M(RB_{\Psi}(\psi_{0}\,|\,D)\geq
1\,|\,\psi_{\ast})$ increases as dist$(\psi_{\ast},\psi_{0})$ increases so
(\ref{biasfor}) is an appropriate measure of bias in favor of $H_{0}$. The
choice of the prior can be used somewhat to control bias but typically a prior
that makes one bias lower just results in making the other bias higher. It is
established in Evans (2105) that, under quite general circumstances, both
biases converge to 0 as the amount of data increases. So bias can be
controlled by design a priori.

The model needs to be checked against the data. For if the data $d$ lies in
the tails of every distribution in the model, then this suggests model
failure. There are a wide variety of approaches to assessing this and these
are not reviewed here. One general comment is that at this time there do not
seem to exist general methodologies for modifying a model when model failure
is encountered.

The prior can also be checked for conflict with the data. A conflict means
that the observed data is in the tails of all those distributions in the model
where the prior primarily places its mass. For a minimal sufficient statistic
$T$ for the model, Evans and Moshonov (2006) used the tail probability
\begin{equation}
M_{T}(m_{T}(t)\leq m_{T}(t)) \label{conflict}%
\end{equation}
to assess prior-data conflict where (\ref{conflict}) small indicates
prior-data conflict. In Evans and Jang (2011a) it is shown that, under general
circumstances, (\ref{conflict}) converges to $\Pi(\pi(\theta)\leq\pi
(\theta_{true}))$ as the amount of data increases. There are a variety of
refinements of (\ref{conflict}) that allow for looking at particular
components of a prior to isolate where a problem with the prior may be. In
Evans and Jang (2011b) a method is developed for replacing a prior when a
prior-data conflict has been detected. This does not mean simply replacing a
prior by one that is more diffuse, however, as is demonstrated in Section 4.1.

\section{Binary-valued Response Regression Models}

The following example, based on real data, is used to illustrate each aspect
of the approach to statistical reasoning recommended here.\smallskip

\noindent\textbf{Example 2. }\textit{Bioassay experiment.}

Table 1 gives the results of exposing animals to various levels in g/ml of a
dosage of a toxin, where $x_{2}$ is the log-dosage and the number of deaths is
recorded at each dosage, see Racine, Grieve, Fluhler and Smith (1986). The
dosages range from $e^{-0.86}=0.423$ to $e^{0.73}=2.075$ g/ml. The logistic
regression model $p(x_{1},x_{2})=G(\beta_{1}+\beta_{2}x_{2})$ is considered
for this data, so $x_{1}\equiv1,G(z)=e^{z}/(1+e^{z}),(\beta_{1},\beta_{2})\in
R^{2}$ and $p(1,x_{2})$ is the probability of death at dosage $x_{2}.$ The
counts $T=(t_{1},t_{2},t_{3},t_{4})$ at the dosages comprise a minimal
sufficient statistic for this problem with observed value $(0,1,3,5).$ The
conditional distribution of $T$ given $(\beta_{1},\beta_{2})$ is a product of binomials.

In Al-Labadai, Baskurt and Evans (2017) a goodness-of-fit test based on this
data was applied for this model using a uniform prior on the space $[0,1]^{4}$
of all probabilities. Relative belief was used to assess the hypothesis that
the model is correct and overwhelming evidence in favor of this model was
obtained and so model correctness is assumed here. One goal is the estimation
of $(\beta_{1},\beta_{2})$ and another is the assessment of the hypothesis
$H_{0}:\beta_{2}=0.$ Acceptance of $H_{0}$ implies that there is no
relationship between the response and the predictor.%

%TCIMACRO{\TeXButton{B}{\begin{table}[tbp] \centering}}%
%BeginExpansion
\begin{table}[tbp] \centering
%EndExpansion%
\begin{tabular}
[c]{|rcc|}\hline
$x_{2}$ & No. of animals & No. of deaths\\\hline
$-0.86$ & 5 & 0\\
$-0.30$ & 5 & 1\\
$-0.05$ & 5 & 3\\
$0.73$ & 5 & 5\\\hline
\end{tabular}
\caption{Data in Example 1.}%
%TCIMACRO{\TeXButton{E}{\end{table}}}%
%BeginExpansion
\end{table}%
%EndExpansion

\subsection{Eliciting the Prior}

\begin{figure}
[ptb]
\begin{center}
\includegraphics[
height=2.303in,
width=2.3142in
]%
{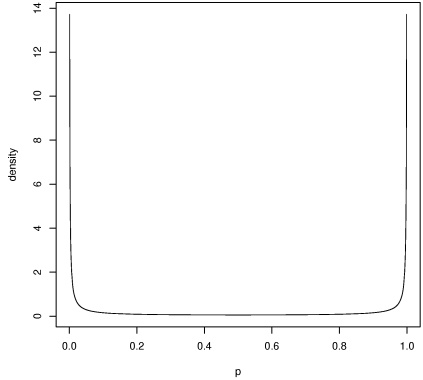}%
\caption{Prior density of of $p(1,x)=G(\beta_{1}+\beta_{2}x),$ $G$ is the
standard logistic cdf, $\beta_{1},\beta_{2}\sim N(0,20^{2})$ and $x=1.$}%
\label{figlogistic1}%
\end{center}
\end{figure}

Elicitation of a prior can be difficult when the interpretation of the
parameters is unclear. For example, with the model (\ref{eq1}) it is not clear
what the $\beta_{i}$ represent in contrast to linear models where they
represent either location parameters or rates of change with respect to
predictors. This leads to attempts to put default priors on these quantities
and there are problems with this approach. For example, suppose
$p(1,x)=G(\beta_{1}+\beta_{2}x),$ where $G$ is the standard logistic cdf and
the prior is given by the $\beta_{i}$ being i.i.d. $N(0,\sigma_{0}^{2})$ where
$\sigma_{0}^{2}$ is chosen large to reflect little information about these
values. In Figure \ref{figlogistic1} we have plotted the prior this induces on
$p(1,1)$ when $\sigma=20.$ This reflects the fact that as $\sigma$ grows all
the prior probability for $p(1,x)$ piles up at 0 and 1 and so this is clearly
a poor choice and it is certainly not noninformative.%
%TCIMACRO{\FRAME{ftbpFU}{2.3142in}{2.303in}{0pt}{\Qcb{Prior density of of
%$p(1,x)=G(\beta_{1}+\beta_{2}x),$ $G$ is the standard logistic cdf, $\beta
%_{1},\beta_{2}\sim N(0,20^{2})$ and $x=1.$}}{\Qlb{figlogistic1}}%
%{unif.eps}{\special{ language "Scientific Word";  type "GRAPHIC";
%width 2.3142in;  height 2.303in;  depth 0pt;  original-width 6.9998in;
%original-height 6.9669in;  cropleft "0";  croptop "1";  cropright "1";
%cropbottom "0";  filename 'unif.eps';file-properties "XNPEU";}} }%
%BeginExpansion
%EndExpansion

The strange behavior of diffuse normal priors has been noted by others.
Bedrick, Christensen and Johnson (1996, 1997), based on Tsutukawa and Lin
(1986),\ make the recommendation that priors should instead be placed on the
$p(\mathbf{x}_{i}),$ as these are parameters for which there is typically
prior information. Their recommendation is that, $k$ of the $\mathbf{x}_{i}$
values be selected and then beta$(\alpha_{1i},\alpha_{2i})$ priors be placed
on the corresponding $p(\mathbf{x}_{i})$ via eliciting prior quantiles$.$ This
results in more sensible priors but depends on the choice of the observed
predictors and it is unclear what kind of priors this induces on the
$\beta_{i}.$

Following Bedrick, Christensen and Johnson (1996, 1997) priors here are
elicited for the probabilities but the approach is different. First, it is not
required that the elicitation be carried out at observed values of the
predictors.\ Rather, it is supposed that there is a set of linearly
independent predictor vectors $\mathbf{w}_{1},\ldots,\mathbf{w}_{k}$ where
bounds can be placed on the probabilities in the sense that $l(\mathbf{w}%
_{i})\leq p(\mathbf{w}_{i})\leq u(\mathbf{w}_{i})$ for $i=1,\ldots,k$ with
virtual certainty. By virtual certainty it is meant that for prior probability
measure $\Pi$, then
\begin{equation}
\Pi(l(\mathbf{w}_{i})\leq p(\mathbf{w}_{i})\leq u(\mathbf{w}_{i})\text{ for
}i=1,\ldots,k)\geq\gamma, \label{ineq1}%
\end{equation}
where $\gamma$ is chosen to be close to 1. For example, $\gamma=0.99$
certainly seems satisfactory for many applications but a higher or lower
standard can be chosen. The motivation for this is that typically information
will be available for the probabilities such as it is known that
$p(\mathbf{w}_{i})$ is very small (or very large) or almost certainly that
$p(\mathbf{w}_{i})$ is in some specific range. Of course, for some of the
$\mathbf{w}_{i}$ virtually nothing may be known about $p(\mathbf{w}_{i})$ and
in that case taking $[l(\mathbf{w}_{i}),u(\mathbf{w}_{i})]=[0,1]$ is
appropriate. One implication of this is that when the choice is made
$[l(\mathbf{w}_{i}),u(\mathbf{w}_{i})]=[0,1]$ for every $i,$ then the
elicitation procedure should lead to a $\Pi$ that is at least approximately
uniform on the probabilities. The approach to elicitation, via stating bounds
on parameter values that hold with virtual certainty, has been successfully
employed in Cao, Evans and Guttman (2015) to determine a prior for the
multivariate normal model, and in Evans, Guttman and Li (2017) to determine a
prior for the multinomial model.

Another reason for allowing the elicitation procedure to be independent of the
observed $\mathbf{x}_{i}$ is that prior beliefs about $p(\mathbf{x}_{i})$ may
apply equally well about $p(\mathbf{x}_{j})$ for some $j$ simply because
$\mathbf{x}_{i}$ and $\mathbf{x}_{j}$ are close and then it seems that the
correlation between the beliefs should be part of the prior. Modelling such
correlations is harder and hopefully can be avoided by choosing the
$\mathbf{w}_{i}$ carefully. For example, requiring the $\mathbf{w}_{i}$ to be
mutually orthogonal seems like an appropriate way of achieving independence in
many contexts.

The second way in which our approach differs from previous developments is
that $\Pi$ is restricted to the family of multivariate normal priors on
$\mathbf{\beta}$ as this allows us to see directly how (\ref{ineq1})
translates into information about $\mathbf{\beta}.$ For note that
(\ref{ineq1}) is equivalent to
\begin{align}
&  \Pi(G^{-1}(l(\mathbf{w}_{i}))\leq G^{-1}(p(\mathbf{w}_{i}))\leq
G^{-1}(u(\mathbf{w}_{i}))\text{ for }i=1,\ldots,k)\nonumber\\
&  =\Pi(G^{-1}(l(\mathbf{w}_{i}))\leq\mathbf{w}_{i}^{\prime}\mathbf{\beta}\leq
G^{-1}(u(\mathbf{w}_{i}))\text{ for }i=1,\ldots,k)\nonumber\\
&  =\Pi(G^{-1}(l(\mathbf{W}))\leq\mathbf{W\beta}\leq G^{-1}(u(\mathbf{W}%
)))\geq\gamma\label{ineq2}%
\end{align}
where $\mathbf{W=(w}_{1}\ldots\mathbf{w}_{k})^{\prime}\in R^{k\times
k},l(\mathbf{W})=\mathbf{(}l\mathbf{(w}_{1}),\ldots,l(\mathbf{w}_{k}%
))^{\prime}\in R^{k},u(\mathbf{W})=\mathbf{(}u\mathbf{(w}_{1}),\ldots
,u(\mathbf{w}_{k}))^{\prime}\in R^{k}.$ So, if $\mathbf{W\beta\sim
}N\mathbf{(\mu}_{0},\Sigma_{0}),$ then $\mathbf{\beta\sim}N_{k}\mathbf{(W}%
^{-1}\mathbf{\mu}_{0},\mathbf{W}^{-1}\Sigma_{0}$\newline$(\mathbf{W}%
^{-1})^{\prime})$ and it is clear what this says about $\mathbf{\beta}.$

The task then is to determine $\mathbf{(\mu}_{0},\Sigma_{0})$ so that
(\ref{ineq2}) is satisfied. A natural choice for $\mathbf{\mu}_{0}$ is to put
$\mathbf{\mu}_{0}=G^{-1}(c(\mathbf{W}))$\ where $c(\mathbf{W})=(l(\mathbf{W}%
)+u(\mathbf{W}))/2$ \ is the centroid of the $k$-cell $[l(\mathbf{W}%
),u(\mathbf{W})].$ For example, when $[l(\mathbf{W}),u(\mathbf{W}%
)]=[0,1]^{k},$ then $c(\mathbf{W})=\mathbf{1}_{k}/2,$ where $\mathbf{1}_{k}$
is the $k$-dimensional vectors of ones, which implies $\mathbf{\mu}%
_{0}=\mathbf{0}.$ Other choices for $\mathbf{\mu}_{0}$ can be made if there
are good reasons for this.

Given that the $\mathbf{w}_{i}$ have been chosen so that prior beliefs about
the probabilities $p(\mathbf{w}_{i})$ are independent, this implies that the
coordinates of $\mathbf{W\beta}$ are independent and so $\Sigma_{0}%
=\,$diag$(\sigma_{1}^{2},\ldots,\sigma_{k}^{2})$ for some choice of the prior
variances $\sigma_{i}^{2}.$ There are, however, typically many choices
satisfying (\ref{ineq2}). For example, taking $\sigma_{i}^{2}=0$ for all $i$
achieves this but clearly this choice does not reflect what is actually known
about the probabilities. As might be expected, the choice of the $\sigma
_{i}^{2}$ is critical and dependent on $G.$ Furthermore, as Figure 1
demonstrates, an injudicious choice results in absurdities.

Since $G^{-1}(u(\mathbf{w}_{i}))-\mu_{0i}>0$ and $G^{-1}(l(\mathbf{w}%
_{i}))-\mu_{0i}<0,$ both these values are infinite iff $[l(\mathbf{w}%
_{i}),u(\mathbf{w}_{i})]=[0,1]$ and so no information is being introduced via
the prior. In such a case a uniform$[0,1]$ prior on the probability results
and the appropriate normal distribution is determined by approximating the
distribution function $G$ by a normal cdf, see Examples 3, 4 and 5. Suppose
then that at least one of $G^{-1}(u(\mathbf{w}_{i}))$ and $G^{-1}%
(l(\mathbf{w}_{i}))$ is finite and so $\sigma_{i}$ satisfies
\begin{equation}
\Phi\left(  \frac{G^{-1}(u(\mathbf{w}_{i}))-\mu_{0i}}{\sigma_{i}}\right)
-\Phi\left(  \frac{G^{-1}(l(\mathbf{w}_{i}))-\mu_{0i}}{\sigma_{i}}\right)
\geq\gamma^{1/k}, \label{marginal}%
\end{equation}
as then independence ensures that (\ref{ineq2}) is satisfied. When both
$G^{-1}(u(\mathbf{w}_{i}))$ and $G^{-1}(l(\mathbf{w}_{i}))$ are finite, the
left side of (\ref{marginal}) has the value 1 when $\sigma_{i}=0,$ is strictly
decreasing to the value 0 as $\sigma_{i}\rightarrow\infty$ and so there are
always values of $\sigma_{i}\geq0$ satisfying (\ref{marginal}). When both
$G^{-1}(u(\mathbf{w}_{i}))$ and $G^{-1}(l(\mathbf{w}_{i}))$ are finite there
is a unique largest solution to (\ref{marginal}), which is the preferred
solution as it best represents the prior information, and it is easily found
numerically by bisection. If $u(\mathbf{w}_{i})=1$ and $l(\mathbf{w}_{i}%
)\in(0,1),$ then $\sigma_{i}=(G^{-1}(l(\mathbf{w}_{i}))-\mu_{0i})/\Phi
^{-1}(1-\gamma^{1/k})$ is the solution provided $\gamma>(1/2)^{k}$ which is a
very weak requirement as recall that $\gamma$ represents virtual certainty. If
$u(\mathbf{w}_{i})\in(0,1)$ and $l(\mathbf{w}_{i})=0,$ then $\sigma
_{i}=(G^{-1}(u(\mathbf{w}_{i}))-\mu_{0i})/\Phi^{-1}(\gamma^{1/k})$ is the
solution again provided $\gamma>(1/2)^{k}.$

The following examples consider the situation $l(\mathbf{w}_{i}%
))=0,u(\mathbf{w}_{i})=1.$ In this case $\mu_{0i}=0$ and $G^{-1}%
(p(\mathbf{w}_{i}))$ will be distributed with cdf $G$ when $p(\mathbf{w}%
_{i})\sim U(0,1).$ Generally this leads to a need to approximate $G$ by\ a
normal cdf to obtain a normal prior although no approximation is required in
Example 3.\smallskip

\noindent\textbf{Example 3.} \textit{Probit regression.}

Here $G=\Phi$ and so $G^{-1}(p(\mathbf{w}_{i}))\sim N(0,1)$ when
$p(\mathbf{w}_{i})\sim U(0,1).$ As such $\sigma_{i}=1$ and the standard normal
distribution on $G^{-1}(p(\mathbf{w}_{i}))$ corresponds to no information
about $p(\mathbf{w}_{i})$. When there is no information about any of the
$p(\mathbf{w}_{i}),$ then $\mathbf{\beta\sim}N_{k}\mathbf{(0},\mathbf{W}%
^{-1}(\mathbf{W}^{-1})^{\prime})$ which equals the $N_{k}\mathbf{(0}%
,\mathbf{I})$ distribution whenever $\mathbf{W}$ is an orthogonal matrix. In
general, however a lack of information about the probabilities leads to a
prior on $\mathbf{\beta}$ that is dependent on $\mathbf{W,}$ namely, dependent
on the values of predictor variables corresponding to the
probabilities.\smallskip

\noindent\textbf{Example 4.} \textit{Logistic regression.}

\begin{figure}
[tb]
\begin{center}
\includegraphics[
height=2.1309in,
width=2.3168in
]%
{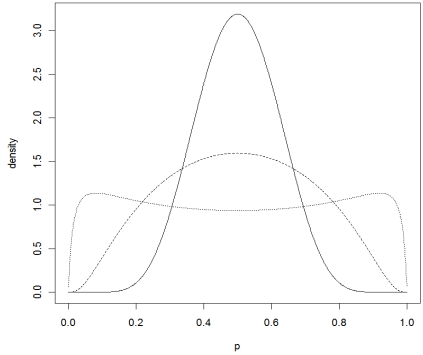}%
\caption{Plots of the density of $p=e^{z}/(1+e^{z})$ when $z\sim N(0,d^{2})$
and $\lambda=0.5$ (--), $\lambda=1.0$ (- -), and $\lambda=-1.702$ (...).
\ \ \ \ \ \ \ \ \ \ \ \ \ \ \ \ \ \ \ \ \ \ \ \ \ \mbox{  }}%
\label{figlogistic2}%
\end{center}
\end{figure}

In this case $G$ is the standard logistic cdf and so $\mathbf{w}_{i}^{\prime
}\mathbf{\beta=}G^{-1}(p(\mathbf{w}_{i}))$ is distributed standard logistic
when $p(\mathbf{w}_{i})\sim U(0,1).$ A well-known $N(0,\lambda^{2})$
approximation to the standard logistic distribution, as discussed in Camilli
(1994), leads to normal priors that are much easier to work with. The optimal
choice of $\lambda,$ in the sense that it minimizes $\max_{x\in R^{1}}%
|\Phi(x/\lambda)-e^{x}/(1+e^{x})|$ is given by $\lambda=1.702$ and this leads
to a maximum difference less than $0.009.$ Clearly this error will generally
be irrelevant when considering priors for the probabilities in a logistic
regression problem. So when $\mathbf{w}_{i}^{\prime}\mathbf{\beta}\sim
N(0,1.702^{2})$ then $p(\mathbf{w}_{i})$ is approximately distributed $U(0,1)$
with the same maximum error. Figure \ref{figlogistic2} contains plots of the
density of $p=e^{z}/(1+e^{z})$ when $z\sim N(0,\lambda^{2})$ for various
choices of $\lambda$ and it is indeed approximately uniform when
$\lambda=1.702.$ Using normal probabilities rather than logistic probabilities
leads to relatively small differences, so it seems reasonable to use a normal
prior on $\mathbf{\beta}$ in a logistic regression\textbf{.}%
%TCIMACRO{\FRAME{ftbFU}{2.3168in}{2.1309in}{0pt}{\Qcb{Plots of the density of
%$p=e^{z}/(1+e^{z})$ when $z\sim N(0,d^{2})$ and $\lambda=0.5$ (--),
%$\lambda=1.0$ (- -), and $\lambda=-1.702$ (...).
%\ \ \ \ \ \ \ \ \ \ \ \ \ \ \ \ \ \ \ \ \ \ \ \ \ \mbox{  }}}%
%{\Qlb{figlogistic2}}{ow6pte00.eps}{\special{ language "Scientific Word";
%type "GRAPHIC";  display "USEDEF";  valid_file "F";  width 2.3168in;
%height 2.1309in;  depth 0pt;  original-width 11.6542in;
%original-height 11.636in;  cropleft "0";  croptop "1";  cropright "1";
%cropbottom "0";  filename 'OW6PTE00.eps';file-properties "XNPEU";}} }%
%BeginExpansion
%EndExpansion
\smallskip

\noindent\textbf{Example 5.} $\mathit{t}$\textit{ regression.}

Suppose that $G$ is taken to be the cdf of $t$ with $\upsilon$ degrees of
freedom. Table 2 presents the optimal choice of $\lambda$ for a $N(0,\lambda
^{2})$ approximation to the $t(\upsilon)$ cdf together with the maximum error.
There does not appear to be much difference in using $t_{\upsilon}$
probabilities instead of normal ones unless $\upsilon$ is quite low.\smallskip%

%TCIMACRO{\TeXButton{B}{\begin{table}[tbp] \centering}}%
%BeginExpansion
\begin{table}[tbp] \centering
%EndExpansion%
\begin{tabular}
[c]{|c|c|c|c|c|c|c|}\hline
$\upsilon$ & $30$ & $20$ & $10$ & $5$ & $2$ & $1$\\\hline
$\lambda$ & \multicolumn{1}{|l|}{$1.022$} & \multicolumn{1}{|l|}{$1.034$} &
\multicolumn{1}{|l|}{$1.069$} & \multicolumn{1}{|l|}{$1.144$} &
\multicolumn{1}{|l|}{$1.407$} & \multicolumn{1}{|l|}{$1.980$}\\
max error & \multicolumn{1}{|l|}{$0.002$} & \multicolumn{1}{|l|}{$0.003$} &
\multicolumn{1}{|l|}{$0.006$} & \multicolumn{1}{|l|}{$0.013$} &
\multicolumn{1}{|l|}{$0.031$} & \multicolumn{1}{|l|}{$0.058$}\\\hline
\end{tabular}
\caption{Optimal choice of a $N(0,\lambda^2)$ distribution to approximate a $t(\nu)$ distribution.}\label{TableKey copy(1)}%
%TCIMACRO{\TeXButton{E}{\end{table}}}%
%BeginExpansion
\end{table}%
%EndExpansion

Consider now an application of the elicitation algorithm.\smallskip

\noindent\textbf{Example 6. }\textit{Bioassay experiment (Example 2
continued).}

In this example $k=2.$ To determine the prior it is necessary to choose
$\mathbf{W=(w}_{1}$ $\mathbf{w}_{2}\mathbf{)}\in R^{2\times2}$ and
$[l(\mathbf{W}),u(\mathbf{W})].$ The authors are not experts in bioassay but,
given the range of dosages applied in the experiment, it is reasonable to
suppose that an expert might be willing to put bounds on the probabilities
that hold with prior probability $\gamma=0.99$ when $x_{2}=-0.50$ and
$x_{2}=0.50$ leading to
\[
\mathbf{W=}\left(
\begin{array}
[c]{cr}%
1 & -1/2\\
1 & 1/2
\end{array}
\right)  .
\]

\begin{figure}
[ptb]
\begin{center}
\includegraphics[
height=2.7345in,
width=2.7345in
]%
{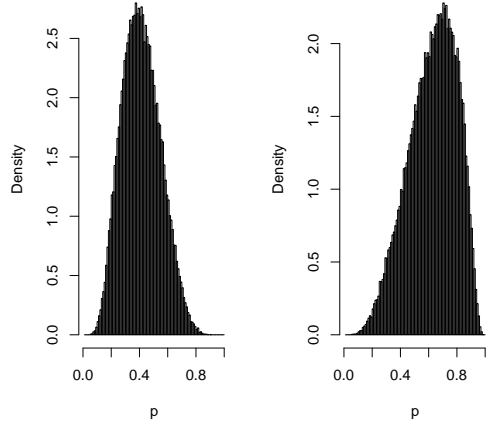}%
\caption{Density histograms of $p(1,-0.8)\,$(left) and $p(1,0.8)$ (right)
based on a sample of $10^{5}$ from the elicited prior in Example 6.}%
\label{priorhists}%
\end{center}
\end{figure}

Let us suppose that an expert believes with virtual certainty that the true
probabilities lie in the intervals $[0.15,0.75],$ when $x_{2}=-0.50,$ and in
$[0.25,0.95],$ when $x_{2}=0.50.$ So the centroid of the $2$-cell
$[0.15,0.75]\times\lbrack0.25,0.95]$ is given by $(0.45,0.60)$ and since
$G^{-1}(p)=\log(p/(1-p))$ for logistic regression, this implies $\mu
_{0}=(G^{-1}(0.45),G^{-1}(0.60))=(-0.20,0.41).$ Also, $[G^{-1}(0.15),G^{-1}%
(0.75]=[-1.735,1.099]$ and $[G^{-1}(0.25),G^{-1}(0.95)]=[-1.099,2.944]$ so,
using (\ref{marginal}), the largest values of $\sigma_{1}$ and $\sigma_{2}$
satisfying $\Phi\left(  1.299/\sigma_{1}\right)  -\Phi\left(  -1.535/\sigma
_{1}\right)  \geq\left(  0.99\right)  ^{1/2}$ and $\Phi\left(  2.534/\sigma
_{2}\right)  -\Phi\left(  -1.509/\sigma_{2}\right)  \geq\left(  0.99\right)
^{1/2}$ are given by $\sigma_{1}=0.490$ and $\sigma_{2}=0.580.$ Therefore the
prior on $\beta$ is
\begin{align}
&  \beta\sim N_{2}(\mathbf{W}^{-1}(-0.20,0.41)^{\prime},\mathbf{W}%
^{-1}\text{diag}(0.490^{2},0.580^{2})(\mathbf{W}^{=1})^{\prime})\nonumber\\
&  =N_{2}\left(  \left(
\begin{array}
[c]{c}%
0.105\\
0.610
\end{array}
\right)  ,\left(
\begin{array}
[c]{cc}%
0.144 & 0.048\\
0.048 & 0.577
\end{array}
\right)  \right)  . \label{prior}%
\end{align}
Figure \ref{priorhists} contains histograms of large samples from the priors
on two extreme probabilities. The shape of the prior is similar for other
values of $x_{2}.$%
%TCIMACRO{\FRAME{ftbpFU}{2.7345in}{2.7345in}{0pt}{\Qcb{Density histograms of
%$p(1,-0.8)\,$(left) and $p(1,0.8)$ (right) based on a sample of $10^{5}$ from
%the elicited prior in Example 6.}}{\Qlb{priorhists}}{fighist.eps}%
%width 2.7345in;  height 2.7345in;  depth 0pt;  original-width 6.9998in;
%original-height 6.9998in;  cropleft "0";  croptop "1";  cropright "1";
%cropbottom "0";  filename 'fighist.eps';file-properties "XNPEU";}} }%
%BeginExpansion
%EndExpansion

\subsection{Measuring the Bias in a Prior}

Consider applying the approach discussed in Section 3.2 to measuring bias in
the prior derived in Section 4.1 for the bioassay example.\smallskip

\noindent\textbf{Example 7. }\textit{Bioassay experiment (Example 2
continued).}

Consider whether or not there is bias induced by the prior in Example 6 with
respect to the hypothesis $H_{0}:\beta_{2}=0.$ It is necessary to compute
$M_{T}(RB_{2}(0\,|\,T)\leq1\,|\,\beta_{2}=0),$ to measure bias against, and
$\sup_{\beta_{2}\in\{-\delta,\delta\}}M_{T}(RB_{2}$\newline$(0\,|\,T)\geq
1\,|\,\beta_{2}),$ to measure bias in favor, where $RB_{2}(\cdot\,|\,T)$ is
the relative belief ratio function for $\beta_{2}$ based on data $T$ and
$\delta>0$ is such that if $|\beta_{2}|<\delta,$ then practically speaking
$H_{0}$ is considered true. To determine $\delta,$ the more general problem of
what changes in both $\beta_{1}$ and $\beta_{2}$ are deemed irrelevant is
considered. Given the settings used in this experiment, it seems reasonable to
consider $x_{2}$ as restricted to the interval $[-1,1].$ Then, whenever
$|\beta_{1}-\beta_{1}^{\prime}|<\delta$ and $|\beta_{2}-\beta_{2}^{\prime
}|<\delta,$ the difference in log odds satisfies $|\beta_{1}+\beta_{2}%
x_{2}-(\beta_{1}^{\prime}+\beta_{2}^{\prime}x_{2})|\leq2\delta$ which implies
that ratio of the odds lies in $(e^{-2\delta},e^{2\delta})$ which for small
$\delta$ is approximately equal to $(1-2\delta,1+2\delta).$ This in turn
implies that the difference in the probabilities is less than $2\delta.$ In
this example we take $\delta=0.01.$

Now $RB(\beta_{1},\beta_{2}\,|\,T)=\{%
%TCIMACRO{\tprod _{i=1}^{4}}%
%BeginExpansion
{\textstyle\prod_{i=1}^{4}}
%EndExpansion
\binom{5}{t_{i}}p^{t_{i}}(1,x_{2i})(1-p(1,x_{2i}))^{5-t_{i}}\}/m_{T}(T)$
where
\[
m_{T}(T)=\int_{-\infty}^{\infty}\int_{-\infty}^{\infty}\{%
%TCIMACRO{\tprod _{i=1}^{4}}%
%BeginExpansion
{\textstyle\prod_{i=1}^{4}}
%EndExpansion
\binom{5}{t_{i}}p^{t_{i}}(1,x_{2i})(1-p(1,x_{2i}))^{5-t_{i}}\}\pi
(\beta)\,d\beta.
\]
The relative belief ratio for $\beta_{2}$ is $RB_{2}(\beta_{2}\,|\,T)=\int
_{-\infty}^{\infty}RB(\beta_{1},\beta_{2}\,|\,T)\pi_{1}(\beta_{1}%
\,|\,\beta_{2})\,d\beta_{1}$\newline$=m_{T}(T\,|\,\beta_{2})/m_{T}(T)$ where
$\pi_{1}(\cdot\,|\,\beta_{2})$ is the conditional prior density of $\beta_{1}$
given $\beta_{2}\ $which (\ref{prior}) implies is the $N(0.105+0.083(\beta
_{2}-0.610),0.140)$ distribution.

For given $T=(t_{1},t_{2},t_{3},t_{4}),$ the numerator and denominator in
$RB_{2}(0\,|\,T)$ can be estimated via simulation but to calculate the biases
we need to do this for many $T.$ Consider the calculation of $M_{T}%
(RB_{2}(0\,|\,T)\leq1\,|\,\beta_{2}=0)$ via the following algorithm and note
that there are only $6^{4}=1296$ values of $(t_{1},t_{2},t_{3},t_{4}%
)\in\{0,1,\ldots,5\}^{4}.$\smallskip

\noindent\textbf{Algorithm}

\noindent(i) simultaneously estimate the values $m_{T}(t_{1},t_{2},t_{3}%
,t_{4})$ for each $(t_{1},t_{2},t_{3},t_{4})$ via a large sample from
(\ref{prior}) and store these in a table,\newline\noindent(ii) simultaneously
estimate the values $m_{T}(t_{1},t_{2},t_{3},t_{4}\,|\,\beta_{2}=0)$ for each
$(t_{1},t_{2},$\newline$t_{3},t_{4})$ via a large sample from $\pi_{1}%
(\cdot\,|\,0)$ and store these in a table,\newline\noindent(iii) using the
values in these two tables estimate $RB_{2}(0\,|\,T)$ for all values of $T$
and then estimate $M_{T}(RB_{2}(0\,|\,T)\leq1\,|\,\beta_{2}=0)$ by summing the
$m_{T}(t_{1},t_{2},t_{3},t_{4}\,|\,$\newline$\beta_{2}=0)$ for those
$(t_{1},t_{2},t_{3},t_{4})$ for which $RB_{2}(0\,|\,T)\leq1.$\smallskip

\noindent The bias in favor can be estimated at $\pm\delta$ in exactly the
same way but in step (ii) replace $\pi_{1}(\cdot\,|\,0)$ by $\pi_{1}%
(\cdot\,|\,-\delta)$ and by $\pi_{1}(\cdot\,|\,\delta)$. These computations
were carried out and resulted in the bias against equaling $0.22$ and the bias
in favor equaling $0.77$ at $-\delta$ and $0.78$ at $\delta$. So there is some
bias against $H_{0}$ with this prior but there is appreciable bias in favor of
$H_{0},$ at least when interest is in detecting deviations of size
$\delta=0.01$. For $\beta_{2}=5,$ however, the bias in favor of $H_{0}$ is
$0.006,$ so there is in reality no bias in favor for large values of this
parameter. One could contemplate modifying the prior to reduce the bias in
favor at $\delta=0.01,$ but typically this just results in trading bias in
favor with bias against. The real cure for excessive bias of either variety,
is to collect more data.\smallskip

In general problems the approach to the computations used here will not be
feasible and so alternative methods are required. In certain examples some
aspects of the computations can be done exactly but, in general,
approximations such as those discussed in Nott et al. (2016) will be necessary.

\subsection{Checking and Modifying a Prior}

Consider now checking the prior derived in Section 4.1 for the bioassay
example.\smallskip

\noindent\textbf{Example 8. }\textit{Bioassay experiment (Example 2
continued).}

The tail probability for checking the prior is given by
\begin{equation}
M_{T}(m_{T}(t_{1},t_{2},t_{3},t_{4})\leq m_{T}(0,1,3,5)). \label{priorchk}%
\end{equation}
As part of the algorithm discussed in Section 4.2, the values of $m_{T}%
(t_{1},t_{2},t_{3},t_{4})$ have been estimated and the proportion of values of
$m_{T}(t_{1},t_{2},t_{3},t_{4})$ that satisfy the inequality gives the
estimate of (\ref{priorchk}). In this example (\ref{priorchk}) equals $0.41$
so there is no prior-data conflict.

If prior-data conflict exists, the methods discussed in Evans and Jang (2011a)
are available to obtain a more weakly informative prior. In this case it is
necessary to be careful as it has been shown in Section 4.1 that simply
increasing the variance of the prior will not necessarily accomplish this. On
the other hand there is the satisfying result that the $N_{2}(\mathbf{0}%
,1.702^{2}I_{2})$ prior$,$ where $I_{2}$ is the identity matrix, will avoid
prior-data conflict, so modifying the elicited prior to be closer to this
prior is the appropriate thing to do when a conflict exists.

\subsection{Inferences}

Now consider estimation and hypothesis assessment for the bioassay
example.\smallskip

\noindent\textbf{Example 9. }\textit{Bioassay experiment (Example 2
continued).}

Consider first the assessment of the hypothesis $H_{0}:\beta_{2}=0.$ From the
algorithm the quantity $RB_{2}(0\,|\,(0,1,3,5))$ is available and this
indicates whether there is evidence in favor of or against $H_{0}.$ In this
case $RB_{2}(0\,|\,(0,1,3,5))=0.021$ so there is evidence against $H_{0}$. A
calculation described below gives the value $0.001$ for the strength, so it
seems there is strong evidence against $H_{0}.$

To obtain the joint relative belief estimate of $(\beta_{1},\beta_{2})$ it is
necessary to maximize $RB(\beta_{1},\beta_{2}\,|\,T)$ as a function of
$(\beta_{1},\beta_{2})$ which is the same as the MLE. The plausibility region
for this estimate is then $\{(\beta_{1},\beta_{2}):RB(\beta_{1},\beta
_{2}\,|\,T)>1\}$ and the size and posterior content of this set provide a
measure of the accuracy with which the coordinates of $\beta$ can be
simultaneously known. But it is worth noting that the $i$-th coordinate of
this joint estimate is not necessarily the value that has the greatest
evidence in its favor, rather this is obtained by maximizing $RB_{i}(\beta
_{i}\,|\,T)$ as a function of $\beta_{i}$ with plausibility region
$\{\beta_{i}:RB_{i}(\beta_{i}\,|\,D)>1\}.$ So the evidence approach dictates
that $\beta_{i}$ be estimated by maximizing $RB_{i}(\beta_{i}\,|\,T).$ In
problems where components of a multidimensional parameter are related by some
constraint, then it is clearly necessary to estimate the components
simultaneously, but that is not the case here.

The value of $RB_{i}(\beta_{i}\,|\,T)$ needs to be estimated and, since this
cannot be done for every value of $\beta_{i},$ its value is estimated on a
finite grid. For this let $[L_{i},U_{i}]$ be the effective prior support for
$\beta_{i}$, say containing $0.995$ of the probability, and form the grid
$G_{i}=\{L_{i},L_{i}+\delta,L_{i}+2\delta,\ldots,U_{i}-\delta,U_{i}\}.$ For
each $\beta_{1}\in G_{1}$ estimate $m_{T}((0,1,3,5)\,|\,\beta_{1})$ using a
large sample from $\pi_{2}(\cdot\,|\,\beta_{1})$ which gives $RB_{1}(\beta
_{1}\,|\,(0,1,3,5))=m_{T}((0,1,3,5)\,|\,\beta_{1})/m_{T}((0,1,3,5)).$ It is
then easy to obtain the relative belief estimate $\beta_{1}(0,1,3,5)$ and
plausibility region $\{\beta_{1}:RB_{1}(\beta_{1}\,|\,(0,1,3,5))>1\}.$ The
true relative belief estimate will differ from this estimate by at most
$\delta$ but this difference has been deemed irrelevant. A similar procedure
is carried out for $\beta_{2}$ but now sampling from $\pi_{1}(\cdot
\,|\,\beta_{2})$ to estimate $m_{T}((0,1,3,5)\,|\,\beta_{2}).$ The posterior
density for $\beta_{i}$ satisfies $\pi_{i}(\beta_{i}\,|\,(0,1,3,5))=RB_{i}%
(\beta_{i}\,|\,(0,1,3,5))\pi_{i}(\beta_{i})$ and since $RB_{i}(\cdot
\,|\,(0,1,3,5))$ has been computed on the grids,\ these values can be used to
approximate the contents of the plausibility regions via an obvious
quadrature. Similarly, the strengths can be estimated and the strength quoted
above equals $%
%TCIMACRO{\tsum _{\beta_{2}\in S\text{ }}}%
%BeginExpansion
{\textstyle\sum_{\beta_{2}\in S\text{ }}}
%EndExpansion
\pi_{2}(\beta_{2}\,|\,(0,1,3,5))\delta$ where $S=G_{2}\cap\{\beta_{2}%
:RB_{2}(\beta_{2}\,|\,(0,1,3,5))\leq RB_{2}(0\,|\,(0,1,3,5))\}.$

\begin{figure}
[ptb]
\begin{center}
\includegraphics[
height=2.4007in,
width=2.4007in
]%
{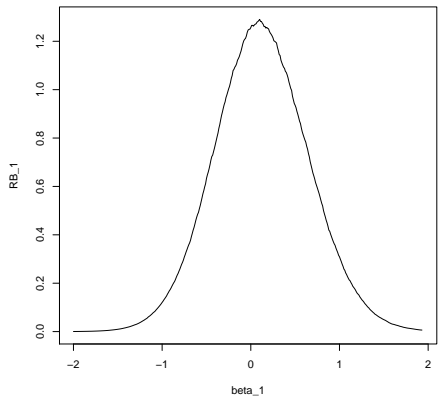}%
\caption{A plot of $RB_{1}(\cdot\,|\,(0,1,3,5))$ over the effective support of
the prior in Example 9.}%
\label{rbbeta1}%
\end{center}
\end{figure}

Implementing this for $\beta_{1},$ the estimate $\beta_{1}(0,1,3,5)=0.11$ was
obtained with plausibility region $[-0.21,0.49]$ having posterior content
$0.35.$ So the range of plausible values for $\beta_{1}$ is not large but
there is not a high belief that the true value is in this interval. Figure
\ref{rbbeta1} is a plot of $RB_{1}(\cdot\,|\,(0,1,3,5)).$%
%TCIMACRO{\FRAME{ftbpFU}{2.4007in}{2.4007in}{0pt}{\Qcb{A plot of $RB_{1}%
%(\cdot\,|\,(0,1,3,5))$ over the effective support of the prior in Example 9.}%
%}{\Qlb{rbbeta1}}{rbbeta1.eps}{\special{ language "Scientific Word";
%type "GRAPHIC";  maintain-aspect-ratio TRUE;  display "USEDEF";
%valid_file "F";  width 2.4007in;  height 2.4007in;  depth 0pt;
%original-width 6.9998in;  original-height 6.9998in;  cropleft "0";
%filename 'rbbeta1.eps';file-properties "XNPEU";}} }%
%BeginExpansion
%EndExpansion

\begin{figure}
[ptb]
\begin{center}
\includegraphics[
height=2.6602in,
width=2.6602in
]%
{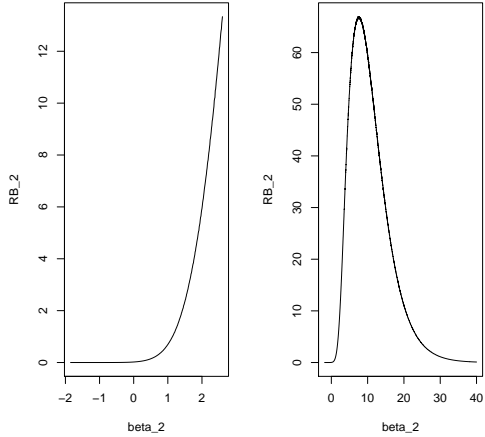}%
\caption{Plot of $RB_{2}(\cdot\,|\,(0,1,3,5)$ over the effective support of
the prior (left panel) and over a full range of possible values (right panel)
in Example 9.}%
\label{rbfig}%
\end{center}
\end{figure}

An interesting phenomenon occurs when considering the estimation of $\beta
_{2}$. In Figure \ref{rbfig} the left panel plots $RB_{2}(\cdot
\,|\,(0,1,3,5))$ over the effective support of the marginal prior $\pi_{2}$
for $\beta_{2}.$ From this it is clear that the relative belief estimate of
$\beta_{2}$ lies outside this range. Recall, however, that the chosen prior
passed the check for prior-data conflict. The check for prior-data conflict
only tells us, however, that the observed data is consistent with at least
some of the probabilities determined by where the prior places its mass. The
right panel of Figure \ref{rbfig} is a plot of $RB_{2}(\cdot\,|\,(0,1,3,5))$
over a much wider range. Note too that there is an important robustness
property as shown in Al-Labadi and Evans (2017) for $RB_{2}(\cdot
\,|\,(0,1,3,5))$ as it is only weakly dependent on $\pi_{2}.$ In this case
$\pi_{2}$ does not place mass where it appears it should but there is not
enough data to detect the conflict. The relative belief estimate of $\beta
_{2}$ is $\beta_{2}(0,1,3,5)=7.31$ and the plausibility region for $\beta_{2}$
is $[1.14,30.48]$ with posterior content $0.83.$ As such, there is a great
deal of uncertainty concerning the true value of $\beta_{2}.$\smallskip%
%TCIMACRO{\FRAME{ftbpFU}{2.6602in}{2.6602in}{0pt}{\Qcb{Plot of $RB_{2}%
%(\cdot\,|\,(0,1,3,5)$ over the effective support of the prior (left panel) and
%over a full range of possible values (right panel) in Example 9.}}%
%{\Qlb{rbfig}}{rbfig.eps}{\special{ language "Scientific Word";
%type "GRAPHIC";  maintain-aspect-ratio TRUE;  display "USEDEF";
%valid_file "F";  width 2.6602in;  height 2.6602in;  depth 0pt;
%original-width 6.9998in;  original-height 6.9998in;  cropleft "0";
%filename 'rbfig.eps';file-properties "XNPEU";}} }%
%BeginExpansion
%EndExpansion

As long as it is possible to sample from the posterior for a 1-dimensional
parameter, then the computations necessary for the inferences for such a
parameter are feasible. As such, the Gibbs sampling algorithm of Albert and
Chib (1993) is particularly relevant although it is not needed in Example 9.
The harder computations are those involving the various prior predictives but
these do not need to be highly accurate as even one decimal place will
indicate whether there is bias or prior-data conflict.

\section{Conclusions}

Criteria for a satisfactory theory of statistical reasoning have been
developed. Perhaps more should be required, but it seems that those stated are
necessary. A particular approach has been outlined that satisfies these
criteria. An example has shown that this approach can resolve
anomalies/paradoxes that arise via a commonly used methodology. Many other
such instances of resolving inferential difficulties as well as results
establishing optimal performance, have been documented in Evans (2015). An
application of the approach to the problem of binary-valued response
regression has been carried out and it has been shown to lead to a number of
novel insights into such problems.

\section{References}

\noindent Albert, J. H. and Chib, S. (1993) Bayesian analysis of binary and
polychotomous response data. J. of the American Statistical Association 88
(422), 669-679.\vspace{2pt}

\noindent Al-Labadi, L. and Evans, M. (2017) Optimal robustness results for
some Bayesian procedures and the relationship to prior-data conflict. Bayesian
Analysis 12, 3, 702-728.\vspace{2pt}

\noindent Al-Labadi, L., Baskurt, Z and Evans, M. (2017) Goodness of fit for
the logistic regression model using relative belief. J. of Statistical
Distributions and Applications, 4:17.\vspace{2pt}

\noindent Baskurt, Z. and Evans, M. (2013) Hypothesis assessment and
inequalities for Bayes factors and relative belief ratios. Bayesian Analysis,
8, 3, 569-590.\vspace{2pt}

\noindent Bedrick, E. J., Christensen, R., and Johnson, W. (1996) A new
perspective on priors for generalized linear models. J. of the American
Statistical Association, 91, 436, 1450-1460.\vspace{2pt}

\noindent Bedrick, E. J., Christensen, R., and Johnson, W. (1997) Bayesian
binomial regression: predicting survival at a trauma center. American
Statistician, 51, 3, 211-218.\vspace{2pt}

\noindent Cao, Y., Evans, M. and Guttman, I. (2015) Bayesian factor analysis
via concentration. Current Trends in Bayesian Methodology with Applications,
edited by S. K. Upadhyay, U. Singh, D. K. Dey and A. Loganathan, 181-201, CRC
Press.\smallskip

\noindent Camilli, G. (1994) Origin of the scaling constant d = 1.7 in item
response theory. J. of Educational and Behavioral Statistics, 19, 3,
293-295.\vspace{2pt}

\noindent Diaconis, P. and Skyrms, B. (2018) Ten Great Ideas About Chance.
Princeton University Press.\vspace{2pt}

\noindent Evans, M. (2015) Measuring Statistical Evidence Using Relative
Belief. Chapman and Hall/CRC.\vspace{2pt}

\noindent Evans, M., Guttman, I. and Li, P. (2017) Prior elicitation,
assessment and inference with a Dirichlet prior. Entropy 2017, 19(10),
564.\smallskip

\noindent Evans, M. and Jang, G-H. (2011a) A limit result for the prior
predictive applied to checking for prior-data conflict. Statistics and
Probability Letters, 81, 1034-1038.\vspace{2pt}

\noindent Evans, M. and Jang, G-H. (2011b) Weak informativity and the
information in one prior relative to another. Statistical Science, 26, 3,
423-439.\vspace{2pt}

\noindent Evans, M. and Moshonov, H. (2006) Checking for prior-data conflict.
Bayesian Analysis, 1, 4, 893-914.\vspace{2pt}

\noindent Kneale, W. and Kneale, M. (1962) The Development of Logic. Clarendon
Pr.\vspace{2pt}

\noindent Nott, D., Drovandi, C., Mengersen, K. and Evans, M. (2016)
Approximation of Bayesian predictive p-values with regression ABC. To appear
in Bayesian Analysis.\vspace{2pt}

\noindent Racine, A., Grieve, A. P., Fluhler, H. and Smith, A. F. M. (1986)
Bayesian methods in practice: experiences in the pharmaceutical industry (with
discussion). J. of Applied Statistics, 35, 93-150.\vspace{2pt}

\noindent Salmon, W. (1973) Confirmation. Scientific American, 228, 5,
75-81.\vspace{2pt}

\noindent Tsutukawa, R. K. and Lin, H. Y. (1986) Bayesian estimation of item
response curves. Psychometrika, 51, 251-267.
\end{document}